%% file: main.tex
\documentclass[12pt,a4paper,leqno]{article}

\usepackage{latexsym}
\usepackage{amssymb,exscale}
\usepackage[centertags]{amsmath}
\usepackage{amsthm}

\usepackage[dvips]{hyperref}

\input{preamble}

\begin{document}

\title{Automorphisms of the type $ I\!I_1 $ Arveson system of Warren's noise}

\author{Boris Tsirelson}

\date{}
\maketitle

\stepcounter{footnote}
\footnotetext{%
 This research was supported by \textsc{the israel science foundation}
 (grant No.~683/05).}

\begin{abstract}
Motions of the plane (shifts and rotations) correspond to automorphisms
of the type $ I $ Arveson system of white noise. I prove that automorphisms
corresponding to rotations cannot be extended to the type $ II $ Arveson
system of Warren's noise.
\end{abstract}

\section*{Introduction}
\input{intro}

\section[Automorphisms of the type $ I_1 $ Arveson system of the white
  noise]{\raggedright Automorphisms of the type $ I_1 $ Arveson system of the
  white noise} 
\label{sect1}
\input{sect1}

\section[An inequality related to CCR]{\raggedright An inequality related to
 CCR}
\label{sect2}
\input{sect2}

\section[The first superchaos of the type $ I\!I_1 $ system]{\raggedright The
 first superchaos of the type $ I\!I_1 $ system}
\label{sect3}
\input{sect3}

\section[Main result]{\raggedright Main result}
\label{sect4}
\input{sect4}

\bigskip
\filbreak
{
\small
\begin{sc}
\parindent=0pt\baselineskip=12pt
\parbox{4in}{
Boris Tsirelson\\
School of Mathematics\\
Tel Aviv University\\
Tel Aviv 69978, Israel
\smallskip
\par\quad\href{mailto:tsirel@post.tau.ac.il}{\tt
 mailto:tsirel@post.tau.ac.il}
\par\quad\href{http://www.tau.ac.il/~tsirel/}{\tt
 http://www.tau.ac.il/\textasciitilde tsirel/}
}

\end{sc}
}
\filbreak

\end{document}

%% file: preamble.tex
\numberwithin{equation}{section}

\swapnumbers
\theoremstyle{definition}
\newtheorem{theorem}[equation]{Theorem}
\newtheorem{lemma}[equation]{Lemma}

\renewcommand{\phi}{\varphi}

\newcommand{\I}{{\rm i}}
\newcommand{\D}{\mathrm{d}}
\newcommand{\E}{\mathrm{e}}

\renewcommand{\(}{\bigl(}
\renewcommand{\)}{\bigr)\vphantom{)}}

\newcommand{\ip}[2]{\langle#1,#2\rangle}

\renewcommand{\equiv}{\;\;\;\Longleftrightarrow\;\;\;}
\newcommand{\thalf}{{\textstyle\frac12}}

\newcommand{\sgn}{\operatorname{sgn}}
\newcommand{\white}{\text{\textup{white}}}
\newcommand{\shift}{\text{\textup{shift}}}
\newcommand{\rotat}{\text{\textup{rotat}}}

\newcommand{\Exp}{\operatorname{Exp}}
\newcommand{\LocMin}{\operatorname{LocMin}}
\newcommand{\Finite}{\operatorname{Finite}}

\newcommand{\One}{\mathbf1}

\newcommand{\eps}{\varepsilon}
\newcommand{\si}{\sigma}
\newcommand{\om}{\omega}
\newcommand{\Om}{\Omega}
\newcommand{\de}{\delta}
\newcommand{\al}{\alpha}
\newcommand{\Ec}{\mathcal E}
\newcommand{\A}{\mathcal A}
\newcommand{\B}{\mathcal B}
\newcommand{\cC}{\mathcal C}
\newcommand{\la}{\lambda}

\newcommand{\Ex}{\mathbb E\,}
\newcommand{\R}{\mathbb R}
\newcommand{\C}{\mathbb C}
\newcommand{\Z}{\mathbb Z}

\newcommand{\sif}{$\sigma$\nobreakdash-field}

\newcommand{\almost}[1]{$#1$\nobreakdash-\hspace{0pt}almost}

%% file: intro.tex
\smallskip

\parbox{4.3cm}{\small\textit{%
\dots we lack information about the gauge groups of non-type I examples.}\\ 
\mbox{}\hfill W.~Arveson \cite[Sect.2.10]{Ar}
}\hfill
\parbox{8.8cm}{\small\textit{%
This is a noise richer than white noise: in addition to the increments of a
Brownian motion $ B $ it carries a countable collection of independent
Bernoulli random variables which are attached to the local minima of $ B $.}\\
\mbox{}\hfill J.~Warren \cite[the end]{Wa}
}

\medskip
\smallskip

The simplest example of a type $ II_1 $ Arveson system emerges from the
simplest example of a nonclassical noise, Warren's noise of splitting (see
\cite{Wa} and \cite{Ts04}, sections 2c, 2e, 4d, 6g). The classical part of the
noise, the white noise, corresponds to the classical part of the Arveson
system, a type $ I_1 $ system. The group of automorphisms of the type $ I_1 $
system, described by Arveson, is basically the group of motions of the plane
(shifts and rotations). It is easy to extend to the nonclassical system the
action of the shifts. It is also easy to see that the rotation by $ \pi $
cannot be extended. However, what happens to other rotations (say, by $ 2\pi/3
$)? It is shown here that only the trivial (by $ 0 $) rotation can be
extended.

%% file: sect1.tex
This section summarizes Arveson's theory of automorphisms of type $ I $
systems, specialized to the white noise, that is, the standard Brownian motion
in $ \R $. The latter is the Gaussian process $ (B_t)_{t\in[0,\infty)} $ with
$ \Ex B_t = 0 $ and $ \Ex (B_s B_t) = s $ for $ 0 \le s \le t < \infty $. Let
$ \Om_t $ be the probability space of the Brownian paths on the time interval
$ [0,t] $. Then
\[
\Om_{s+t} = \Om_s \times \Om_t
\]
up to the natural isomorphism of probability spaces. Thus,
\[
L_2 (\Om_{s+t}) = L_2 (\Om_s) \otimes L_2(\Om_t)
\]
up to the natural isomorphism of Hilbert spaces (just a unitary operator). It
means that these $ L_2(\Om_t) $ form an Arveson system (as defined in
\cite[Sect.~3.1, Def.~3.1.1]{Ar}).

The space $ L_2(\Om_t) $ may be thought of as the exponential of $ L_2(0,t)
$ (see \cite[Sect.~2.1]{Ar}, especially (2.7)). To this end we define a map $
\Exp : L_2(0,t) \to L_2(\Om_t) $ by 
\[
\Exp (f) = \exp \bigg( \int_0^t f(s) \, \D B_s - \frac12 \int_0^t f^2(s) \, \D
s \bigg)
\]
(here $ \exp x $ means the usual $ \E^x $, not to be confused with $ \Exp $),
using the well-known unitary operator $ f \mapsto \int f \, \D B $ from $
L_2(0,t) $ onto a subspace of $ L_2(\Om_t) $ (in fact, the subspace of all
measurable linear functionals of $ B $). Linear combinations of these $
\Exp(f) $ are dense in $ L_2(\Om_t) $, and
\[
\ip{ \Exp(f) }{ \Exp(g) } = \exp \ip f g \, ,
\]
since $ \Ex \( \exp ( \int f \, \D B - \frac12 \int f^2 \, \D s ) \,
\overline{ \exp ( \int g \, \D B - \frac12 \int g^2 \, \D s ) } \, \) = \exp
\( -\frac12 \int (f^2+\overline{g}^2) \, \D s \) \Ex \exp \( \int
(f+\overline{g}) \, \D B \) $, and $ \Ex \exp \( \int (f+\overline{g}) \, \D B
\) = \exp \( \frac12 \int (f+\overline{g})^2 \, \D s \) = \exp \( \frac12 \int
(f^2+\overline{g}^2) \, \D s \) \exp \( \int f \overline{g} \, \D s \) $. It
means that the Arveson system is (isomorphic to) the exponential product
system of rank $ 1 $ \cite[Sect.~3.1, Def.~3.1.6]{Ar}, therefore it is of
type $ I_1 $ (in other words, of type $ I $ and index $ 1 $).

The units of the Arveson system (see \cite[Sect.~3.6]{Ar}, especially (3.21))
are of the form $ u^{(a,\zeta)} $ for $ a, \zeta \in \C $; here
\begin{equation}\label{1.1}
u^{(a,\zeta)} (t) = \E^{at} \Exp ( \zeta \cdot \chi_{(0,t)} ) = \exp \(
\zeta B_t - \thalf \zeta^2 t + at \) \, .
\end{equation}

\begin{sloppypar}
Isomorphisms between Arveson systems are defined in \cite[Sect.~3.1,
Def.~3.1.3]{Ar}; isomorphisms to itself are called automorphisms. An
automorphism $ \theta $ consists of unitary operators $ \theta_t : L_2(\Om_t)
\to L_2(\Om_t) $ such that $ \theta_{s+t} = \theta_s \otimes \theta_t $. The
automorphisms of the type $ I_1 $ system (see \cite[Sect.~3.8]{Ar}, especially
(3.37) and the proof of Th.~3.8.4) are of the form $ \theta^{(\la,\xi,U)} $
for $ \la \in \R $, $ \xi \in \C $, $ U \in \C $, $ |U| = 1 $; they act on the
units as follows:
\begin{gather*}
\theta^{(\la,\xi,U)} u^{(a,\zeta)} = u^{(a',\zeta')} \, , \\
\zeta' = U\zeta + \xi \, , \quad a' = a + \I\la - \thalf |\xi|^2 - U \zeta
\overline\xi \, .
\end{gather*}
See \cite[Sect.~3.8]{Ar} for the composition formula (for two automorphisms)
in terms of $ \la,\xi,U $. 
\end{sloppypar}

Parameters $ a $ and $ \la $ are trivial in the sense that
\[
u^{(a,\zeta)} (t) = \E^{at} u^{(\zeta)} (t) \, , \quad \theta^{(\la,\xi,U)}_t
= \E^{\I\la t} \theta^{(\xi,U)}_t
\]
where $ u^{(\zeta)} = u^{(0,\zeta)} $ and $ \theta^{(\xi,U)} =
\theta^{(0,\xi,U)} $. Accordingly,
\begin{equation}
\theta^{(\xi,U)}_t u^{(\zeta)} (t) = \exp \( -\thalf |\xi|^2 t - U
\zeta \overline\xi t \) u^{(U\zeta+\xi)} (t) \, .
\end{equation}
Denoting for convenience
\[
\theta^{\shift(\xi)} = \theta^{(\xi,1)} \, , \quad \theta^{\rotat(U)} =
\theta^{(0,U)}
\]
we have
\begin{gather}
\theta^{\rotat(U)} \theta^{\rotat(V)} = \theta^{\rotat(UV)} \, , \quad
 \theta^{\rotat(U)} u^{(\zeta)} = u^{(U\zeta)} \, , \\
\theta^{\rotat(U)} \theta^{\shift(\xi)} \( \theta^{\rotat(U)} \)^{-1} =
 \theta^{\shift(U\xi)} \, , \label{1.4} \\
\theta^{\shift(\xi)}_t u^{(\zeta)} (t) = \exp \( -\thalf |\xi|^2 t -
\zeta \overline\xi t \) u^{(\zeta+\xi)} (t) \, .
\end{gather}
For $ \la,\mu \in \R $,
\begin{gather}
\theta^{\shift(\I\la)} \theta^{\shift(\I\mu)} = \theta^{\shift(\I(\la+\mu))}
 \, , \quad \theta^{\shift(\la)} \theta^{\shift(\mu)} =
 \theta^{\shift(\la+\mu)} \, , \\
\theta^{\shift(\I\la)}_t u^{(\zeta)} (t) = \exp \( -\thalf \la^2 t + \I
 \la \zeta t \) u^{(\zeta+\I\la)} (t) \, , \label{1.7} \\
\theta^{\shift(\la)}_t u^{(\zeta)} (t) = \exp \( -\thalf \la^2 t - \la
 \zeta t \) u^{(\zeta+\la)} (t) \, ,
\end{gather}
which leads to canonical commutation relations (CCR) \cite[Remark 3.8.2]{Ar}
\begin{equation}\label{1.9}
\theta^{\shift(\I\la)}_t \theta^{\shift(\mu)}_t = \E^{2\I\la\mu t}
\theta^{\shift(\mu)}_t \theta^{\shift(\I\la)}_t \, .
\end{equation}
Combining \eqref{1.1} and \eqref{1.7} we get $ \theta^{\shift(\I\la)}_t
\exp(\zeta B_t) = \exp (\I\la B_t) \exp (\zeta B_t) $, thus identifying the
automorphism $ \theta^{\shift(\I\la)} $ with the automorphism formed by
multiplication operators,
\begin{equation}\label{1.10}
\theta^{\shift(\I\la)}_t f = \exp (\I\la B_t) f \quad \text{for } f \in
L_2(\Om_t) \, .
\end{equation}

%% file: sect2.tex
Let $ P,Q $ be selfadjoint operators on a separable Hilbert space. The
canonical commutation relations $ [P,Q]=-\I $ will be treated as an
abbreviation of the Weyl relations
\[
\forall \la,\mu \in \R \quad \E^{\I\la P} \E^{\I\mu Q} = \E^{\I\la\mu}
\E^{\I\mu Q} \E^{\I\la P} \, .
\]
If $ [P,Q]=-\I $ then $ P+Q $ is well-defined and $ [Q,-(P+Q)]=-\I $. Thus, we
may speak about three operators $ P,Q,R $ such that $ P+Q+R = 0 $ and $
[P,Q]=-\I $, $ [Q,R]=-\I $, $ [R,P]=-\I $ (these three relations being in fact
mutually equivalent).

\begin{theorem}\label{2.1}
Let selfadjoint operators $ P,Q,R $ be such that $ P+Q+R = 0 $ and $
[P,Q] = [Q,R] = [R,P] =-\I $. Then
\[
\| \sgn P + \sgn Q + \sgn R \| < 3 \, .
\]
\end{theorem}

(Here `$ \sgn P $' is the discontinuous sign function applied to the operator
$ P $.) The proof is given in \cite[Th.~2.1]{Ts06} for the irreducible
representation of CCR (unique up to unitary equivalence). The general case
follows easily, since every representation decomposes into irreducible
representations (von Neumann's theorem).

Note that $ \| \sgn P + \sgn Q + \sgn R \| $ is an absolute constant (since
all irreducible triples $ P,Q,R $ are mutually unitarily
equivalent). According to a numerical computation \cite[Sect.~1]{Ts06}, the
constant is approximately $ 2.1 $.

Returning to the context of Sect.~\ref{sect1} we introduce the generator $
Q_t $ of the unitary group $ \( \theta^{\shift(\I\la)}_t \)_{\la\in\R} $,
\[
\E^{\I\la Q_t} = \theta^{\shift(\I\la)}_t \, ;
\]
by \eqref{1.10}, it is the multiplication by $ B_t : \Om_t \to \R $,
\begin{equation}\label{2.2}
Q_t f = B_t f \quad \text{for $ f \in L_2(\Om_t) $ such that $ B_t f \in
L_2(\Om_t) $} \, .
\end{equation}
By \eqref{1.4}, the operator $ \theta^{\rotat(U)}_t Q_t \(
\theta^{\rotat(U)}_t \)^{-1} $ is the generator of the unitary group $ \(
\theta^{\shift(\I U\la)}_t \)_{\la\in\R} $. Especially, the operator
\[
P_t = \theta^{\rotat(\I)}_t Q_t \theta^{\rotat(-\I)}_t
\]
satisfies
\[
\E^{-\I\mu P_t} = \theta^{\shift(\mu)}_t \, ,
\]
and we may rewrite \eqref{1.9} as
\[
[ P_t, Q_t ] = -2t\I \, .
\]
More generally, for $ \al \in \R $
\begin{equation}\label{2.25}
\theta^{\rotat(\E^{\I\al})}_t Q_t \theta^{\rotat(\E^{-\I\al})}_t = Q_t \cos\al
+ P_t \sin\al \, .
\end{equation}

\begin{lemma}\label{2.3}
There exists $ \eps > 0 $ such that for every $ t \in (0,\infty) $, $ \al \in
\(\frac\pi2,\pi] $ and $ f \in L_2(\Om_t) $,
\[
\langle \sgn Q_t \rangle_f + \langle \sgn Q_t \rangle_g + \langle \sgn Q_t
\rangle_h \le ( 3 - \eps ) \| f \|^2 \, ;
\]
here $ g = \theta^{\rotat(\E^{\I\al})}_t f $, $ h =
\theta^{\rotat(\E^{-\I\al})}_t f $, and $ \langle A \rangle_f $ stands for $
\ip{Af}{f} $.
\end{lemma}

\begin{proof}
Using the general relations
\[
\langle A \rangle_g = \langle \theta^{\rotat(\E^{-\I\al})}_t A
\theta^{\rotat(\E^{\I\al})}_t \rangle_f
\]
and
\[
\theta^{\rotat(\E^{-\I\al})}_t (\sgn A) \theta^{\rotat(\E^{\I\al})}_t = \sgn
\( \theta^{\rotat(\E^{-\I\al})}_t A \theta^{\rotat(\E^{\I\al})}_t \)
\]
we get by \eqref{2.25}
\begin{gather*}
\langle \sgn Q_t \rangle_g = \langle \sgn ( Q_t \cos\al - P_t \sin\al )
 \rangle_f \, , \\
\langle \sgn Q_t \rangle_h = \langle \sgn ( Q_t \cos\al + P_t \sin\al )
 \rangle_f \, .
\end{gather*}
If $ \al=\pi $ then $ \langle \sgn Q_t \rangle_f + \langle \sgn Q_t \rangle_g
= 0 $, thus, the inequality holds (for $ \eps=2 $). Otherwise,
Theorem \ref{2.1} may be applied to the operators $ P = a Q_t $, $ Q = b ( Q_t
\cos\al - P_t \sin\al ) $, $ R = c ( Q_t \cos\al + P_t \sin\al ) $ provided
that $ a,b,c \in (0,\infty) $ are chosen appropriately (namely, $ a = \sqrt{
(-\cos\al) / (t\sin\al) } $ and $ b=c= 1/\sqrt{-4t\sin\al\cos\al} $). We get
\begin{multline*}
\langle \sgn Q_t \rangle_f + \langle \sgn ( Q_t \cos\al + P_t \sin\al )
 \rangle_f + \langle \sgn ( Q_t \cos\al - P_t \sin\al ) \rangle_f \le \\
\le \| \sgn Q_t + \sgn ( Q_t \cos\al + P_t \sin\al ) + \sgn ( Q_t \cos\al -
 P_t \sin\al ) \| \cdot \| f \|^2 \le \\
\le (3-\eps) \|f\|^2 \, ,
\end{multline*}
where $ 3-\eps $ is the absolute constant given by Theorem \ref{2.1}.
\end{proof}

%% file: sect3.tex
Probability spaces denoted by $ \Om_t $ in sections \ref{sect1}, \ref{sect2}
will be denoted by $ \Om_t^\white $ in sections \ref{sect3},
\ref{sect4}. Similarly, other objects relating to the white noise will be
marked `white', because we turn to Warren's noise of splitting, richer than
the white noise.

A path $ \om_{1} $ of the noise of splitting on the time interval $ [0,1] $
consists of a Brownian path $ \om^\white_{1} \in \Om^\white_{1} $ and a
map $ \eta_{1} : \LocMin (\om^\white_{1}) \to \{-1,+1\} $; here $ \LocMin
(\om^\white_{1}) $ is the set of all local minimizers of the path $
\om^\white_{1} : [0,1] \to \R $, and $ \Om_{1}^\white \subset C_0[0,1] $
is a Borel set of full Wiener measure such that for every $ \om^\white_{1} \in
\Om_{1}^\white $ the set $ \LocMin (\om^\white_{1}) $ is a dense countable
subset of $ (0,1) $, and all the local minima are strict. We may choose a
measurable enumeration of local minimizers, that is, a sequence of measurable
maps $ \tau_1, \tau_2, \dots : \Om_{1}^\white \to (0,1) $ such that 
\[
\LocMin (\om^\white_{1}) = \{ \tau_1(\om^\white_{1}),
\tau_2(\om^\white_{1}), \dots \}
\]
for every $ \om^\white_{1} \in \Om_{1}^\white $, and these $
\tau_k(\om^\white_{1}) $ are pairwise different.

Every measurable enumeration $ (\tau_k)_k $ of the local minimizers on $ (0,1)
$ gives us a one-to-one correspondence
\begin{gather*}
\Om_{1} \leftrightarrow \Om_{1}^\white \times \{-1,+1\}^\infty \, , \\
\om_{1} = ( \om^\white_{1}, \eta_{1} ) \leftrightarrow \(
\om^\white_{1}, \( \eta_{1} ( \tau_1(\om^\white_{1})), \eta_{1} (
\tau_2(\om^\white_{1})), \dots \) \) \, ;
\end{gather*}
here $ \Om_1 $ is the set of all paths of the noise of splitting on the time
interval $ [0,1] $, and $ \{-1,+1\}^\infty $ is the set of all infinite
sequences of $ \pm1 $. We equip $ \{-1,+1\}^\infty $ with the product measure
$ m^\infty $, where $ m $ gives to $ -1 $ and $ +1 $ equal probabilities $
1/2, 1/2 $. Further, we equip $ \Om_{1}^\white \times \{-1,+1\}^\infty $ with
the Wiener measure multiplied by $ m^\infty $. Finally, using the one-to-one
correspondence, we transfer the probability measure (and the underlying \sif)
to $ \Om_{1} $, getting $ P_1 $. The choice of an enumeration $ (\tau_k)_k $
does not matter, since $ m^\infty $ is invariant under permutations.

Probability spaces $ \Om_t = (\Om_t,P_t) $ for $ t \in (0,\infty) $ are
constructed similarly; they satisfy $ \Om_{s+t} = \Om_s \times \Om_t $.

The general form of a function $ f \in L_2 ( \Om_{1}^\white \times
\{-1,+1\}^\infty ) $ is
\begin{gather*}
f \( \om^\white_{1}, (\si_1,\si_2,\dots) \) = \sum_{n=0}^\infty
 \sum_{k_1<\dots<k_n} \si_{k_1}\dots\si_{k_n} f_{k_1,\dots,k_n}
 (\om^\white_{1}) \, , \\
f_{k_1,\dots,k_n} \in L_2(\Om^\white_{1}) \, , \quad \sum_{n=0}^\infty
 \sum_{k_1<\dots<k_n} \| f_{k_1,\dots,k_n} \|^2 = \| f \|^2 < \infty \, ;
\end{gather*}
of course, $ \si_k = \pm1 $. Therefore the general form of $ f \in L_2 (
\Om_{1} ) $ is
\begin{gather}
f(\om_{1}) = \sum_{n=0}^\infty \sum_{k_1<\dots<k_n} \eta_{1} (
 \tau_{k_1} (\om^\white_{1})) \dots \eta_{1} (\tau_{k_n}
 (\om^\white_{1})) f_{k_1,\dots,k_n} (\om^\white_{1}) \, , \label{C.1} \\
f_{k_1,\dots,k_n} \in L_2(\Om^\white_{1}) \, , \quad \sum_{n=0}^\infty
 \sum_{k_1<\dots<k_n} \| f_{k_1,\dots,k_n} \|^2 = \| f \|^2 < \infty \,
 . \notag
\end{gather}
For $ n=0 $ we get the natural embedding $ L_2 ( \Om_{1}^\white ) \subset
L_2 ( \Om_{1} ) $.

The Hilbert spaces $ H_t = L_2 ( \Om_{t} ) $ for $ t \in (0,\infty) $ are an
Arveson system. Its automorphisms $ \theta $ consist of unitary operators $
\theta_t : H_t \to H_t $. The subspace $ H_t^\white = L_2 ( \Om_t^\white ) $
of $ H_t $ is invariant under $ \theta_t $ (for every automorphism $ \theta $) 
since, first, the classical (in other words: type $ I $; completely spatial;
decomposable) part of an Arveson system is invariant under automorphisms, and
second, the classical part of the system $ (H_t)_t $ is the system $
(H_t^\white)_t $ (see \cite[Sections 4d, 6g]{Ts04}). 

The set $ \Finite(0,1) $ of all finite subsets of $ (0,1) $ is a Borel
space. Every bounded Borel function $ \phi : \Finite(0,1) \to \R $ leads to an
operator $ \Ec_\phi : H_1 \to H_1 $, given in terms of \eqref{C.1} by
\begin{multline*}
(\Ec_\phi f) (\om_{1}) = \sum_{n=0}^\infty \sum_{k_1<\dots<k_n} \eta_{1} (
 \tau_{k_1} (\om^\white_{1})) \dots \eta_{1} (\tau_{k_n}
 (\om^\white_{1})) \cdot \\
\cdot f_{k_1,\dots,k_n} (\om^\white_{1}) \phi \( \{ \tau_{k_1}
 (\om^\white_{1}), \dots, \tau_{k_n} (\om^\white_{1}) \} \) \, .
\end{multline*}
Thus, the commutative algebra of all bounded Borel functions on $ \Finite(0,1)
$ acts on $ H_1 $. Its action commutes with automorphisms; this fact is a
special case of a more general statement \cite[Sect.~3]{Ts04a}, but I give a
streamlined proof here.

\begin{lemma}
Operators $ \Ec_\phi $ and  $ \theta_1 $ commute, whenever $ \theta =
(\theta_t)_{t\in(0,\infty)} $ is an automorphism of the Arveson system $
(H_t)_{t\in(0,\infty)} $ and $ \phi : \Finite(0,1) \to \R $ is a bounded Borel
function.
\end{lemma}

\begin{proof}
The orthogonal projection onto the subspace $ H_s \otimes H^\white_{t-s}
\otimes H_{1-t} \subset H_1 $ (for $ 0<s<t<1 $) commutes with $ \theta_1 =
\theta_s \otimes \theta_{t-s} \otimes \theta_{1-t} $ and is of the form $
\Ec_\phi $; namely, $ \phi(C)=1 $ if $ C \cap (s,t) = \emptyset $, otherwise $
\phi(C) = 0 $. Thus, the lemma holds for these special $ \phi $.

The Borel sets of the form $ \{ C \in \Finite(0,1) : C \cap (s,t) = \emptyset
\} $ generate the Borel \sif\ of $ \Finite(0,1) $. Proof: restricting
ourselves to rational $ s, t $ we get a countable collection of Borel sets
separating points of $ \Finite(0,1) $, therefore, generating the Borel
\sif.

It means that the lemma holds for all $ g $ taking on the values $ 0, 1 $
only. The general case follows.
\end{proof}

Thus, a subspace of $ H_1 $ corresponds to every Borel subset of $
\Finite(0,1) $. Especially, the classical part, $ H^\white_{1} $,
corresponds to $ \{ C : |C|=0 \} = \{\emptyset\} $ (just $ n=0 $ in
\eqref{C.1}). The subspace corresponding to $ \{ C : |C|=1 \} = \{ \{t\} :
0<t<1 \} $ (just $ n=1 $ in \eqref{C.1}) is the so-called \emph{first
superchaos space} $ H_1^{(1)} \subset H_1 $;
\begin{gather*}
f \in H_1^{(1)} \equiv f(\om_{1}) = \sum_k \eta_{1} ( \tau_k (\om^\white_{1}))
 f_k (\om^\white_{1}) \, , \\
f_k \in L_2(\Om^\white_{1}) \, , \quad \sum_k \| f_k \|^2 = \| f \|^2 <
 \infty \, .
\end{gather*}
(Subspaces $ H_t^{(1)} \subset H_t $ appear similarly.)
Automorphisms leave $ H_1^{(1)} $ invariant. The commutative algebra of
bounded Borel functions $ \chi : (0,1) \to \R $ acts on $ H_1^{(1)} $,
\[
(\A_\chi f) (\om_{1}) = \sum_k \eta_{1} ( \tau_k (\om^\white_{1})) f_k
(\om^\white_{1}) \chi ( \tau_k (\om^\white_{1}) ) \, ,
\]
and commutes with automorphisms restricted to $ H_1^{(1)} $. For example,
taking $ \chi (\cdot) = 1 $ on $ (0,t) $ and $ \chi (\cdot) = 0 $ on $ (t,1) $
we get the orthogonal projection onto the subspace $ H_t^{(1)} \otimes
H_{1-t}^\white \subset H_1 $.

Operators $ \A_\chi $ commute also with the natural action of the algebra $
L_\infty (\Om_1^\white) $ on $ H_1^{(1)} $,
\[
(\B_\phi f) (\om_{1}) = \sum_k \eta_{1} ( \tau_k (\om^\white_{1})) f_k
(\om^\white_{1}) \phi (\om^\white_{1}) \, .
\]
However, the action $ \phi \mapsto \B_\phi $ does not commute with
automorphisms. We may join the actions $ \A,\B $ into an action $ \cC $ of the
commutative algebra of bounded Borel functions $ \psi : (0,1) \times
\Om_1^\white \to \R $ on $ H_1^{(1)} $,
\[
(\cC_\psi f) (\om_{1}) = \sum_k \eta_{1} ( \tau_k (\om^\white_{1})) f_k
(\om^\white_{1}) \psi ( \tau_k (\om^\white_{1}), \om^\white_{1} ) \, .
\]
In particular, consider the function
\begin{equation}\label{3.3}
\psi (t,\om^\white_1) = \begin{cases}
 \sgn \( B_1(\om^\white_1)-B_{0.5}(\om^\white_1) \) &\text{if $ t<0.5 $}, \\
 0 &\text{otherwise}.
\end{cases}
\end{equation}
(Of course, $ B_t(\om^\white_1) $ is just another notation for $
\om_1^\white(t) $.)
It acts on $ H_1^{(1)} = H_{0.5}^{(1)} \otimes H_{0.5}^\white \oplus
H_{0.5}^\white \otimes H_{0.5}^{(1)} $ as follows (recall \eqref{2.2}):
\begin{equation}\label{3.4}
\begin{aligned}
\cC_\psi &= \One_{0.5} \otimes \sgn Q_{0.5} &\text{on } H_{0.5}^{(1)} \otimes
 H_{0.5}^\white \, , \\
\cC_\psi &= 0 &\text{on } H_{0.5}^\white \otimes H_{0.5}^{(1)} \, .
\end{aligned}
\end{equation}

Each function $ f \in H_1^{(1)} $ leads to a finite positive Borel measure $
\mu_f $ on $ (0,1) \times \Om_1^\white $ such that for every bounded Borel $
\psi $,
\[
\int \psi \, \D\mu_f = \langle \cC_\psi \rangle_f = \sum_k \int_{\Om_1^\white}
| f_k (\om^\white_{1}) |^2 \psi ( \tau_k (\om^\white_{1}), \om^\white_{1} ) \,
P_1^\white (\D\om^\white_{1}) \, .
\]
Clearly, $ \mu_f \( (0,1) \times \Om_1^\white \) = \|f\|^2 $, and $ t \in
\LocMin (\om^\white_{1}) $ for \almost{\mu_f} all pairs $ (t,\om^\white_{1})
$.

%% file: sect4.tex
\begin{theorem}\label{4.1}
For every $ \al \in (0,2\pi) $ the automorphism $ \theta^{\rotat(\E^{\I\al})}
$ of the classical part $ (H_t^\white)_t $ of the Arveson system $ (H_t)_t $
cannot be extended to an automorphism of the whole system.
\end{theorem}

Assume the contrary: the extension $ \theta $ exists for some $ \al \in
(0,2\pi) $. We also assume that $ \al \in \(\frac\pi2,\pi] $ (otherwise we may
use $ n\al $ for an appropriate $ n \in \Z $). As before, $ \langle A
\rangle_f $ stands for $ \ip{Af}{f} $. The operator $ Q_t $ acts on $
H_t^\white $, recall \eqref{2.2}.

\begin{lemma}\label{4.2}
Let $ s,t > 0 $, $ f \in H_s^{(1)} \otimes H_t^\white \subset H_{s+t} $. Then
\[
\langle \One_s \otimes \sgn Q_t \rangle_f + \langle \One_s \otimes \sgn Q_t
\rangle_g + \langle \One_s \otimes \sgn Q_t \rangle_h \le (3-\eps) \| f \|^2
\, ,
\]
where $ \One_s $ is the identity operator on $ H_s $, $ g = \theta_{s+t} f $,
$ h = \theta_{s+t}^{-1} f $ and $ \eps $ is the same as in Lemma \ref{2.3} (a
positive absolute constant).
\end{lemma}

\begin{proof}
We repeat the proof of Lemma \ref{2.3}, taking into account that $
\theta_{s+t}^{-1} (\One_s \otimes Q_t) \theta_{s+t} = \( \theta_s^{-1} \otimes
\theta_t^{\rotat(\E^{-\I\al})} \) \( \One_s \otimes Q_t \) \( \theta_s \otimes
\theta_t^{\rotat(\E^{\I\al})} \) = \One_s \otimes (Q_t \cos\al - P_t \sin\al)
$ on $ H_s^{(1)} \otimes H_t^\white $.
\end{proof}

The following construction is the key to the proof of Theorem \ref{4.1}. For
any $ n \in \{1,2,\dots\} $ and $ \de \in (0,0.5) $ we define Borel functions
$ \psi_{n,\de} : (0,1) \times \Om_1^\white \to [-1,1] $ by
\[
\psi_{n,\de} (t,\om_1^\white) = \sum_{k=1}^n \chi_{n,k}(t) \phi_{n,k,\de}
(\om_1^\white) \, ,
\]
where
\begin{gather*}
\chi_{n,k}(t) = \begin{cases}
  1 &\text{if } t \in [\frac{k-1}{2n},\frac{k}{2n}\), \\
  0 &\text{otherwise};
 \end{cases} \\
\phi_{n,k,\de} (\om_1^\white) = \sgn \( B_{\frac{k}{2n}+\de} (\om_1^\white) -
 B_{\frac{k}{2n}} (\om_1^\white) \) \, .
\end{gather*}
Note that the functions $ \chi_{n,k}(t) \phi_{n,k,\de} (\om_1^\white) $ are
similar to \eqref{3.3} and act similarly to \eqref{3.4}.

\begin{lemma}\label{4.3}
For every $ f \in H_{0.5}^{(1)} \otimes H_{0.5}^\white $,
\[
\liminf_{n\to\infty} \langle \cC_{\psi_{n,\de}} \rangle_f \to \| f \|^2 \quad
  \text{as } \de \to 0+ \, .
\]
\end{lemma}

\begin{proof}
We define Borel sets $ U_\de \subset (0,0.5) \times \Om_1^\white $ by
\[
(t,\om_1^\white) \in U_\de \equiv B_{t+\de}(\om_1^\white) > B_t(\om_1^\white)
\, .
\]
For \almost{\mu_f} all pairs $ (t,\om_1^\white) $ we have $ t \in (0,0.5) \cap
\LocMin(\om_1^\white) $, therefore, $ (t,\om_1^\white) \in U_\de $ for all $
\de $ small enough. It follows that
\[
\mu_f (U_\de) \to \| f \|^2 \quad \text{as } \de \to 0+ \, .
\]
If $ (t,\om_1^\white) \in U_\de $ then (by continuity of Brownian paths) $
\psi_{n,\de} (t,\om_1^\white) = +1 $ for all $ n $ large enough. Therefore
\[
\int_{U_\de} \psi_{n,\de} \, \D\mu_f \to \mu_f (U_\de) \quad \text{as } n \to
\infty \, .
\]
However, $ \langle \cC_{\psi_{n,\de}} \rangle_f = \int \psi_{n,\de} \, \D\mu_f
\ge \int_{U_\de} \psi_{n,\de} \, \D\mu_f - \mu_f ( \complement U_\de ) $
(since $ \psi_{n,\de}(\cdot) \ge -1 $), thus,
\[
\liminf_{n\to\infty} \langle \cC_{\psi_{n,\de}} \rangle_f \ge \mu_f(U_\de) -
\mu_f ( \complement U_\de ) \to \| f \|^2 \quad \text{as } \de \to 0+ \, .
\]
Also, $ \langle \cC_{\psi_{n,\de}} \rangle_f \le \| f \|^2 $ (since $
\psi_{n,\de}(\cdot) \le 1 $).
\end{proof}

\begin{lemma}\label{4.4}
For all $ f \in H_{0.5}^{(1)} \otimes H_{0.5}^\white $ and all $ n,\de $
\[
\langle \cC_{\psi_{n,\de}} \rangle_f + \langle \cC_{\psi_{n,\de}} \rangle_g +
\langle \cC_{\psi_{n,\de}} \rangle_h \le (3-\eps) \|f\|^2 \, ,
\]
where $ g = \theta_1 f $, $ h = \theta_1^{-1} f $, and $ \eps $ is the same as
in Lemma \ref{2.3} (a positive absolute constant).
\end{lemma}

\begin{proof}
Applying Lemma \ref{4.2} (or rather, its evident generalization) to $
\A_{\chi_{n,k}} f $ (in place of $ f $) and taking into account that $
\A_{\chi_{n,k}} g = \theta_1 \A_{\chi_{n,k}} f $, $ \A_{\chi_{n,k}} h =
\theta_1^{-1} \A_{\chi_{n,k}} f $ (since $ \A_{\chi_{n,k}} $ commutes with $
\theta_1 $) we get
\[
\langle \A_{\chi_{n,k}} \B_{\phi_{n,k,\de}} \rangle_f + \langle \A_{\chi_{n,k}}
\B_{\phi_{n,k,\de}} \rangle_g + \langle \A_{\chi_{n,k}} \B_{\phi_{n,k,\de}}
\rangle_h \le (3-\eps) \| \A_{\chi_{n,k}} f \|^2 \, .
\]
We sum up in $ k $ and note that $ \sum_k \A_{\chi_{n,k}} \B_{\phi_{n,k,\de}} =
\cC_{\psi_{n,\de}} $ and $ \sum_k \| \A_{\chi_{n,k}} f \|^2 = \| f \|^2 $.
\end{proof}

Applying Lemma \ref{4.3} to $ f $, $ g = \theta_1 f $ and $ h = \theta_1^{-1}
f $ we get
\[
\liminf_{n\to\infty} \( \langle \cC_{\psi_{n,\de}} \rangle_f + \langle
\cC_{\psi_{n,\de}} \rangle_g + \langle \cC_{\psi_{n,\de}} \rangle_h \) \to 3
\| f \|^2 \quad \text{as } \de \to 0+
\]
in contradiction to Lemma \ref{4.4}, which completes the proof of Theorem
\ref{4.1}.

%% file: main.bbl
\begin{thebibliography}{8.}

{\raggedright
\bibitem{Ar} W. Arveson (2003):
\emph{Noncommutative dynamics and $E$-semigroups,}
Springer, New York.

\bibitem{Ts04} B. Tsirelson (2004):
\emph{Nonclassical stochastic flows and continuous products,}
Probability Surveys \textbf{1}, 173--298.

\bibitem{Ts04a} B. Tsirelson (2004):
\emph{On automorphisms of type II Arveson systems \textup{(}probabilistic
approach\textup{),}} 
arXiv:math.OA/0411062v1.

\bibitem{Ts06} B. Tsirelson (2006):
\emph{How often is the coordinate of a harmonic oscillator positive?},
arXiv:quant-ph/0611147v1.

\bibitem{Wa} J. Warren (1999):
\emph{Splitting: Tanaka's SDE revisited,}
arXiv:math.PR/9911115v1.


}
\end{thebibliography}
